\documentclass[11pt]{article}
\usepackage{graphicx}
\usepackage{amssymb}
\usepackage{amsmath}
\usepackage{amsthm}
\newtheorem{theorem}{Theorem}[section]
\newtheorem{lemma}[theorem]{Lemma}
\newtheorem{proposition}[theorem]{Proposition}
\newtheorem{corollary}[theorem]{Corollary}
\newtheorem{definition}[theorem]{Definition}

\def \C{\mathbb{C}}
\def \P{{\sf I\kern-1.5ptP}}

\newcommand{\Ree}{{\rm Re}}
\newcommand{\Imm}{{\rm Im}}
\begin{document}
\thispagestyle{plain}
\begin{center} \Large
     The asymptotic properties of the spectrum of non symmetrically perturbed
     Jacobi matrix sequences
  \\[20pt]
  \large
   Leonid Golinskii and Stefano Serra-Capizzano
  \footnote{The work of the first author is supported in part by INTAS Research network
  03-51-6637. The work of the second author is supported in part by MIUR grant no. 2002014121.}
   \end{center}
\begin{abstract}
Under the mild trace-norm assumptions, we show that the
eigenvalues of a generic (non Hermitian) complex perturbation of a
Jacobi matrix sequence (not necessarily real) are still
distributed as the real-valued function $2\cos t$ on $[0,\pi]$
which characterizes the nonperturbed case. In this way the real
interval $[-2,2]$ is still a cluster for the asymptotic joint
spectrum and, moreover, $[-2,2]$ still attracts strongly (with
infinite order) the perturbed matrix sequence. The results follow
in a straightforward way from more general facts that we prove
in an asymptotic linear algebra framework and are plainly
generalized to the case of matrix-valued symbols, which arises
when dealing with orthogonal polynomials with asymptotically
periodic recurrence coefficients.
\end{abstract}
\hfill
\\ {\bf AMS Classification (2000): 15A18, 15A12, 47B36, 47B65}

\section{Introduction and preliminary discussion}\label{sec:intro}

Consider the matrix $J_n^0$ of size $n$ defined as
\begin{equation} \label{toep-conv}
J_n^0= \left[ \begin{array}{cccccc}
    0 & 1 & \\
    1 & 0 & 1 & \\
        & 1 &\ddots& \ddots \\
        &    &\ddots & \ddots & 1\\
        &     &      & 1 & 0
     \end{array} \right].
\end{equation}
The former matrix is the Toeplitz matrix $T_n(a)$ generated by
$a(t)=2\cos t$ in the following sense: given a Lebesgue integrable
function $b$ defined on $[-\pi,\pi)$ (and periodically extended on
$\mathbb R$), the matrix $T_n(b)$ has order $n$ and entries
$\left(T_n(b)\right)_{p,q}={\hat b}_{p-q}$,\ $p,q=1,\ldots,n$.
Here ${\hat b}_{j}$ is the $j$-th Fourier coefficient of $b$,
i.e.,
\[
{\hat b}_{j}=\frac1{2\pi}\,\int_{-\pi}^\pi b(t) {\exp}(-ijt)\, dt,
\qquad j \in \mathbb Z,\qquad i^2=-1.
\]
In the specific case (\ref{toep-conv}) the eigenvalues are
explicitly known, and they coincide with the evaluation of $a(t)$
on the uniform grid $j\pi/(n+1)$ on $[0,\pi]$. If $J_n^0$ is
replaced by a more general Jacobi matrix
\begin{equation} \label{trid-gen-conv}
J_n = \left[ \begin{array}{cccccc}
    b_0 & a_{1} & \\
    a_{1} & b_1 & a_{2} & \\
        & a_{2} &\ddots& \ddots \\
        &    &\ddots & \ddots & a_{n-1} \\
        &     &      & a_{n-1} & b_{n-1}
     \end{array} \right],
\end{equation}
where $a_j\in\mathbb R$ tends to $1$ and $b_j\in\mathbb R$ tends to $0$ as
$j\to\infty$, then its eigenvalues are no longer explicitly known,
but they are again an approximation of the evaluation of $a(t)$ on
the same grid. This result can be obtained directly from the GLT
theory (see \cite{ser-glt,glt-vs-fourier}), and more precisely,
$\forall F\in \mathcal C_0(\mathbb C)$ ($\mathcal C_0(\mathbb C)$
is the set of all continuous functions having a bounded support),
we have
\begin{equation} \label{def_asym-bis}
     \lim_{n \to \infty} \frac{1}{n} \sum_{\lambda\in \Sigma_n} F(\lambda)
     = \frac1{\pi}\,\int_0^\pi F(2\cos t)\,dt=
     \frac1{2\pi}\,\int_{-\pi}^\pi F(2\cos t)\,dt.
\end{equation}
Here and in what follows $\Sigma_n$ stands for the collection of
all eigenvalues of $J_n$ counted with their multiplicity, the
function $2\cos t$ is also called the symbol of $\{J_n\}$, and we
write $\{J_n\}\sim_\lambda (2\cos t,[-\pi,\pi])$. In the
orthogonal polynomials community this result is known, often under the
unnecessary condition $a_j>0$, in the form
\[
\lim_{n\to\infty} \frac1n\sum_{j=1}^n
F(x_{j,n})=\frac1{\pi}\int_{-2}^2 \frac{F(x)\,dx}{\sqrt{4-x^2}}\,,
\]
that is the weak*-convergence of the counting measures of the
zeros $\{x_{j,n}\}_{j=1}^n$ of orthonormal polynomials $\{p_n\}$
to the equilibrium measure of the support of the orthogonality
measure (see, e.g., \cite [Section 4.9]{nevai}, \cite [Chapter
2]{stahltot}). Observe that the set of zeros $\{x_{j,n}\}_{j=1}^n$
is exactly the set $\Sigma_n$ considered in the left-hand side of
(\ref{def_asym-bis}).

Let us set up the formal definitions. For any function $F$ defined
on $\mathbb C$ and any matrix $A_n$ of size $n$, with the
eigenvalues $\lambda_j(A_n)$, $j=1,\ldots,n$, the symbol
$\Sigma_{\lambda}(F,A_n)$ stands for the mean
\[
\Sigma_{\lambda}(F,A_n):= {\frac 1n \sum_{j=1}^{n}
F\left(\lambda_j(A_n)\right)}={\frac 1 {n} \sum_{\lambda\in \Sigma_n}
F(\lambda)}.
\]
A generic sequence of matrices $\{A_n\}:=\{A_n\}_n$ ($A_n$ of size $n$)
will be referred to as a {\em matrix sequence}.
\begin{definition}\label{def-distribution}
{\rm A matrix sequence $\{A_n\}$ is {\em distributed $($in the
sense of the eigenvalues$)$ as a measurable function $\theta$},
defined on a set $G\subset\mathbb R^q$ of finite and positive
Lebesgue measure $m(G)$, if $\forall F\in \mathcal C_0(\mathbb
C)$, the following limit relation holds
\begin{equation}\label{distribution:sv-eig}
\lim_{n\rightarrow
\infty}\Sigma_{\lambda}(F,A_n)=\frac1{m(G)}\,\int_G F(\theta(t))\,
dt.
\end{equation}
In this case we write in short $\{A_n\}\sim_{\lambda} (\theta,G)$.
Moreover, two sequences $\{A_n\}$ and $\{B_n\}$ are {\em equally
distributed} if $\forall F\in \mathcal C_0(\mathbb C)$, we have
\begin{equation}\label{distribution:equal}
\lim_{n\rightarrow
\infty}[\Sigma_{\lambda}(F,B_n)-\Sigma_{\lambda}(F,A_n)]=0.
\end{equation}}
\end{definition}
Note that two sequences having the same distribution function are
equally distributed. On the other hand, two equally distributed
sequences do not need to have a distribution function. However, if
one of them has a distribution function then the other necessarily
shares the same distribution: the derivation is immediate from the
definitions (for an example see \cite [Remark 6.1]{taud2}).

Along with the distribution in the sense of eigenvalues
(weak*-convergence) we will study another asymptotic property of
the spectra $\Sigma_n$ called here the {\em clustering}.
\begin{definition}\label{def-cluster}
{\rm A matrix sequence $\{A_n\}$ is {\em properly $($or
strongly$)$ clustered at $s \in \mathbb C$} (in the eigenvalue
sense) if for any $\varepsilon>0$ the number of the eigenvalues of
$A_n$ off the disk
\[
D(s,\varepsilon):=\{z:|z-s|<\varepsilon\}
\]
can be bounded by a pure constant $q_\varepsilon$ possibly
depending on $\varepsilon$, but not on $n$. In other words
\[
q_\varepsilon(n,s):=\#\{\lambda_j(A_n): \lambda_j\notin
D(s,\varepsilon)\}=O(1), \quad n\to\infty.
\]
If every $A_n$ has, at least definitely (that is, for all large
enough $n$), only real eigenvalues, then $s$ is real and the disk
$D(s,\varepsilon)$ reduces to the interval
$(s-\varepsilon,s+\varepsilon)$. Furthermore, $\{A_n\}$ is {\em
properly $($or strongly$)$ clustered at a nonempty closed set $S
\subset \mathbb C$} (in the eigenvalue sense) if for any
$\varepsilon>0$
\begin{equation} \label{2.2}
q_\varepsilon(n,S):=\#\{\lambda_j(A_n): \lambda_j\not\in
D(S,\varepsilon):=\cup_{s\in S} D(s,\varepsilon)\}=O(1), \quad
n\to\infty,
\end{equation}
$D(S,\varepsilon)$ is the $\varepsilon$-neighborhood of $S$, and
if every $A_n$ has, at least definitely, only real eigenvalues,
then $S$ has to be a nonempty closed subset of $\mathbb R$.
Finally, the term ``properly (or strongly)'' is replaced by
``weakly'', if
\[
q_\varepsilon(n,s)=o(n), \qquad
\bigl(q_\varepsilon(n,S)=o(n)\bigr), \quad n\to\infty,
\]
in the case of a point $s$ (a closed set $S$), respectively.}
\end{definition}
It is clear that $\{A_n\}\sim_{\lambda} (\theta,G)$ with
$\theta\equiv s$ a constant function is equivalent to $\{A_n\}$
being weakly clustered at $s \in \mathbb C$ (for more results and
relations among the notions of equal distribution, equal
localization, spectral distribution, spectral clustering etc., see
\cite [Section 4]{taud2}).

We will primarily be interested in the special situation, when
$J_n$ are viewed as $n\times n$ principal blocks of an infinite
Jacobi matrix $J_\infty$ (background), $P_\infty$ is a complex
Jacobi matrix (perturbation), $A_\infty=J_\infty+P_\infty$ and
$A_n=J_n+P_n$ are the $n\times n$ principal blocks of $A_\infty$
(so $A_{n+1}$ is the one step extension of $A_n$). In fact, the
main results hold in much more general setting when no relation
between $A_{n+1}$ and $A_n$ is presumed.

The main conditions we impose on $P_\infty$ are of two types.
\\
(i). $\|P_n\|_1=o(n)$ as $n\to\infty$, where $\|\cdot\|_1$ is the
trace-norm of a matrix (i.e., the sum of its singular values, see
\cite{bhatia}). This condition is equivalent to
\begin{equation}\label{equiv-1}
\lim_{n\to\infty} \frac1n
\sum_{j=1}^n\bigl(|p_{j,j-1}|+|p_{j,j}|+|p_{j,j+1}|\bigr)=0, \quad
P_\infty=\{p_{j,k}\}_{j,k=1}^\infty.
\end{equation}
The latter means the Ces\`aro convergence of the entries of
$P_\infty$ to zero. $P_\infty$ is now called the {\em Ces\`aro
compact Jacobi matrix} (cf. \cite{gha, golin}). \\
(ii). $\|P_n\|_1=O(1)$ as $n\to\infty$, that is,
\begin{equation}\label{equiv-2}
\limsup_{n\to\infty}\sum_{j=1}^n
\bigl(|p_{j,j-1}|+|p_{j,j}|+|p_{j,j+1}|\bigr)<\infty,
\end{equation}
and so $A_\infty$ is the trace class perturbation of $J_\infty$.

We point out that the trace-norm is useful in the theoretical
derivations while the conditions on the entries are easy to check
in practice. Moreover, the equivalence of the trace-norm and
entry-wise conditions in (i) and (ii) is well known (cf.
\cite[Section 2]{kilsimon}). Nevertheless we give the proof in
Appendix for two reasons, i.e., because we get better equivalence
constants, for the sake of completeness, and because the proposed
matrix-theoretic proof is new and elementary.

We proceed as follows. In Section \ref{sec:tools} the relation
between the distribution in the sense of eigenvalues, clustering
and attracting properties of matrix sequences is discussed. Our
main results are stated and proved in Section \ref{sec:results}.
In particular, Theorems~\ref{mirsky-cons} and
\ref{mirsky-cons-bis} allow to study non-Hermitian perturbations
of Hermitian matrix sequences. As a straightforward consequence we
obtain the clustering for zeros of the system of polynomials
satisfying the three-term recurrence relation with complex
coefficients. Finally, in Section \ref{sec:periodic} we examine
the case of block Toeplitz and asymptotically periodic Jacobi
matrices, and in Section \ref{sec:extension} we discuss further
extensions and generalizations.

\section{Clustering and attracting}\label{sec:tools}

Let us recall the notion of the essential range which plays an
important role in the study of asymptotic properties of the
spectrum.
\begin{definition}\label{def-range}
{\rm Given a measurable complex-valued function $\theta$ defined
on a Lebesgue measurable set $G$, the {\em essential range of
$\theta$} is the set $S(\theta)$ of points $s\in \mathbb C$ such
that, for every $\varepsilon>0$, the Lebesgue measure of the set
$\theta^{(-1)}(D(s,\varepsilon)):=\{t\in G:\ \theta(t)\in
D(s,\varepsilon)\}$ is positive. The function $\theta$ is {\em
essentially bounded} if its essential range is bounded. Finally,
if $\theta$ is real-valued, then the essential supremum (infimum)
is defined as the supremum (infimum) of its essential range.}
\end{definition}
$S(\theta)$ is clearly a closed set (its complement is open), and
moreover
$$ S(\theta)=\bigcap \{B-{\rm closed\ set}:
m(\theta^{(-1)}(B))=m(G)\}, $$ where $m(X)$ is the Lebesgue
measure of a set $X$.

In the case of a bounded in the operator norm sequence $\{A_n\}$,
a further mathematical instrument that we need is a way for
relating formula (\ref{distribution:sv-eig}), with $F$ a generic
polynomial, to the same formula in its full extent, i.e., with $F$
being a continuous function. The answer is partly contained in the
Mergelyan Theorem and not completely positive. We need assumptions
on the essential range of the symbol $\theta$ and a priori
assumptions on the clustering properties of the sequence
$\{A_n\}$. The reason is in part due to the barrier given by the
Mergelyan Theorem stating that the closure in the uniform norm of
the polynomials on a compact set $S$ is given by the set of all
continuous functions on $S$ which are holomorphic in its interior,
provided that $\mathbb C\backslash S$ is connected (for the proof
see \cite [Theorem 20.5, pp. 423-427]  {Ru}). Therefore the
polynomial space is able to approximate every continuous function
on $S$ if and only if $S$ has the empty interior and $\mathbb
C\backslash S$ is connected.

\begin{theorem}\label{mergelyan-cons}
Assume that a matrix sequence $\{A_n\}$ is weakly clustered at a
compact set $S\subset\mathbb C$ with the connected complement, and
the spectra $\Sigma_n$ are uniformly bounded, i.e., $|\lambda|<C$,
$\lambda\in\Sigma_n$, for all $n$. Assume further that
$(\ref{distribution:sv-eig})$ holds with $F$ a generic polynomial
of an arbitrary fixed degree, and the essential range of $\theta$
in contained in $S$. Then relation $(\ref{distribution:sv-eig})$
is true for every continuous function $F$ with a bounded support
which is holomorphic in the interior of $S$. Moreover, if the
interior of $S$ is empty, then $\{A_n\}$ is distributed as
$\theta$ on its domain $G$.
\end{theorem}
\noindent
 {\em Proof.} In the argument we follow Tilli (see
\cite{tillicomplex}, the proof of Theorem 3). Take $F$ continuous
over $S$ and holomorphic in its interior. By the Mergelyan
Theorem, for every $\varepsilon>0$, we can find a polynomial $p$
such that $|p(z)-F(z)|\le \varepsilon$ for every $z\in S$. Since
the essential range of $\theta$ is contained in $S$, it is clear
that $|p(\theta(t))-F(\theta(t))|\le \varepsilon$ a.e. in its
domain $G$. Therefore
\begin{equation}\label{rel1}
 \left|\frac1{m(G)}\int_G F(\theta(t))dt-\frac1{m(G)}
 \int_G p(\theta(t))dt\right|
 \le \frac{\varepsilon}{m(G)}\int_G dt=\varepsilon.
\end{equation}
Next, we go over to the left-hand side of
(\ref{distribution:sv-eig}). By the definition of clustering for
any fixed $\varepsilon'>0$ we have
\[
\#\{\lambda\in \Sigma_n,\ |\lambda-z|\ge \varepsilon',\ \forall
z\in S\}= \#\{\lambda\in \Sigma_n,\ \lambda\notin
D(S,\varepsilon')\}=o(n).
\]
Moreover, by the hypothesis of the uniform boundedness of
$\Sigma_n$, $|\lambda|<C$ for every $\lambda\in \Sigma_n$ with a
pure constant $C$ independent of $n$. Therefore, by extending $F$
outside $S$ in such a way that it is continuous with a bounded
support, we have
\[
\left| \frac{1}{n} \sum_{\lambda\in \Sigma_n,\ \lambda\notin
D(S,\varepsilon')} F(\lambda)\right| \le \frac{M}{n}
\#\{\lambda\in \Sigma_n,\ \lambda\notin D(S,\varepsilon')\}=o(1),
\]
\[
\left| \frac{1}{n} \sum_{\lambda\in \Sigma_n,\ \lambda\notin
D(S,\varepsilon')} p(\lambda)\right| \le \frac{M}{n}
\#\{\lambda\in \Sigma_n,\ \lambda\notin D(S,\varepsilon')\}=o(1),
\]
with $M=\max(\|F\|_\infty, \|p\|_\infty)$, and the infinity norms
are taken over $\{z\in \C,\ |z|\le C\}$. Consequently, by setting
$\displaystyle \Delta=\left|\Sigma_{\lambda}(F-p,A_n)\right|$ we
have
\begin{eqnarray*}
\Delta &=& \left| \frac{1}{n} \sum_{\lambda\in \Sigma_n}
(F(\lambda)-p(\lambda))\right| \le \frac{1}{n} \sum_{\lambda\in
\Sigma_n} |F(\lambda)-p(\lambda)| \\
{} &=& \frac{1}{n} \sum_{\lambda\in \Sigma_n,\ \lambda\in
D(S,\varepsilon')} |F(\lambda)-p(\lambda)|+ \frac{1}{n}
\sum_{\lambda\in \Sigma_n,\ \lambda\notin D(S,\varepsilon')}
|F(\lambda)-p(\lambda)|
\\
{} &\le & \frac{1}{n} \sum_{\lambda\in \Sigma_n,\ \lambda\in
D(S,\varepsilon')} |F(\lambda)-p(\lambda)|+o(1)
\\
{} &=& \frac{1}{n} \sum_{\lambda\in \Sigma_n,\ \lambda\in S}
|F(\lambda)-p(\lambda)| +\frac{1}{n} \sum_{\lambda\in \Sigma_n,\
\lambda\in D(S,\varepsilon')\backslash S}
       |F(\lambda)-p(\lambda)|+o(1).
\end{eqnarray*}
For $\lambda\in S$ we use $|F(\lambda)-p(\lambda)|\le
\varepsilon$, and for $\lambda\in D(S,\varepsilon')\backslash S$
we write
\[
|F(\lambda)-p(\lambda)|\le |F(\lambda)-F(\lambda')|+
|F(\lambda')-p(\lambda')|+|p(\lambda')-p(\lambda)|,\ \ \
|\lambda-\lambda'|< \varepsilon', \ \ \lambda'\in S,
\]
so that $|F(\lambda)-p(\lambda)|\le
c_1(\varepsilon')+\varepsilon+c_2(\varepsilon,\varepsilon')\equiv
\theta(\varepsilon,\varepsilon')$ with
\begin{equation}\label{rel2}
\lim_{\varepsilon\rightarrow 0} \lim_{\varepsilon' \rightarrow 0}
\theta(\varepsilon,\varepsilon')=0.
\end{equation}
Hence
\begin{equation} \label{rel3}
\Delta\le \varepsilon+\theta(\varepsilon,\varepsilon')+o(1).
\end{equation}
 Moreover, from the hypothesis of the theorem we have
\begin{equation} \label{asym-pol}
  \lim_{n \to \infty}\Sigma_{\lambda}(p,A_n)
     = \frac1{m(G)}\int_G p(\theta(t))dt.
\end{equation}
Since $\varepsilon$ and $\varepsilon'$ are arbitrary, it is clear
that relations (\ref{rel1})--(\ref{asym-pol}) imply
(\ref{distribution:sv-eig}) to hold for $F$ as well. Finally, when
$S$ has empty interior, we have no restriction on $F$ except for
being continuous with a bounded support, and therefore what we
have proved is equivalent to $\{A_n\}\sim_\lambda (\theta,G)$.
\hfill $\square$

To proceed further, we need a notion which is essential in the
orthogonal polynomials theory.

\begin{definition}\label{def-attractors}
{\rm A matrix sequence $\{A_n\}$ is {\em strongly attracted by $s
\in \mathbb C$} if
\begin{equation} \label{sattr}
\lim_{n\rightarrow \infty}{\rm dist}(s,\Sigma_n)=0,
\end{equation}
where ${\rm dist}(X,Y)$ is the usual Euclidean distance between
two subsets $X$ and $Y$ of the complex plane. Furthermore, let us
order the eigenvalues according to its distance from $s$, i.e.,
\[
|\lambda_1(A_n)-s|\le |\lambda_2(A_n)-s| \le \cdots \le
|\lambda_n(A_n)-s|.
\]
We say that the attraction is of order
$r(s)\in \mathbb N$, $r(s)\ge 1$ is a fixed number, if
\[
\lim_{n\rightarrow \infty}|\lambda_{r(s)}(A_n)-s|=0,\ \ \ \
\liminf_{n\rightarrow \infty}|\lambda_{r(s)+1}(A_n)-s|>0.
\]
The attraction is of order $r(s)=\infty$ if
\[
\lim_{n\rightarrow \infty}|\lambda_{j}(A_n)-s|=0
\]
for every fixed $j$. Finally, the term ``strong or strongly'' is
replaced by ``weak or weakly'' if $\lim$ is replaced by $\liminf$
in (\ref{sattr}).}
\end{definition}

It is not hard to ascertain, that if $\{A_n\}$ is at least weakly
clustered at a point $s$, then $s$ strongly attracts $\{A_n\}$
with infinite order. Indeed, $s$ is an attracting point of finite
order implies
\[
\lim_{n\to\infty}\frac{\#\{\lambda\in \Sigma_n:\ \lambda\notin
D(s,\delta)\}}{n}=1
\]
for some $\delta>0$, that is impossible in the case when $\{A_n\}$
is weakly clustered at $s$. On the other hand, there are sequences
which are strongly attracted by $s$ with infinite order but not
even weakly clustered at $s$.

The notions previously introduced in this section are intimately
related, as emphasized in the following theorem.

\begin{theorem}\label{th-relations-def}
Let $\theta$ be a measurable function defined on $G$ with finite
and positive Lebesgue measure, and $S=S(\theta)$ the essential
range of $\theta$. Let $\{A_n\}$ be a matrix sequence distributed
as $\theta$ in the sense of eigenvalues. Then
\begin{description}
 \item[a)]  $S(\theta)$ is a weak cluster for $\{A_n\}$;
 \item[b)] each point $s\in S(\theta)$ strongly attracts $\Sigma_n$ with
infinite order $r(s)=\infty$.
\end{description}
\end{theorem}
\noindent {\em Proof.} {\bf a)}. Given $\varepsilon>0$, we apply
(\ref{distribution:sv-eig}) with the test function $F_\varepsilon$
of the form
\[
F_\varepsilon(z)=
  \begin{cases}
    1, & \text{for} \ z\in D(S,\varepsilon/2)\cap D(0,1/\varepsilon), \\
    0, & \text{for} \ z\in \mathbb C\backslash \left(D(S,\varepsilon)\cap
        D(0,2/\varepsilon)\right),
  \end{cases} \qquad 0\leq F_\varepsilon\leq 1.
\]
It is clear that
\[
\begin{aligned} \Sigma_\lambda(F_\varepsilon, A_n) &\leq
\frac{\#\{\lambda\in\Sigma_n: \lambda\in
\left(D(S,\varepsilon)\cap
    D(0,2/\varepsilon)\right)\}}{n} \\
    &\leq \frac{\#\{\lambda\in\Sigma_n: \lambda\in
    D(S,\varepsilon)\}}{n}=1-\frac{q_\varepsilon(n,S)}{n}\,,
    \end{aligned}\]
$q_\epsilon(n,S)$ is defined in (\ref{2.2}), and so
\[
\liminf_{n\to\infty} \Sigma_\lambda(F_\varepsilon, A_n)\leq
1-\limsup_{n\to\infty}\frac{q_\varepsilon(n,S)}{n}.
\]
 By the assumption there exists
\[
\lim_{n\to\infty} \Sigma_\lambda(F_\varepsilon,
A_n)=\frac1{m(G)}\int_G F(\theta(t))\,dt \geq
\frac{m\{\theta^{(-1)}\left(D(S,\varepsilon/2)\cap
    D(0,1/\varepsilon)\right)\}}{m(G)}\,. \]
We have
\[
\theta^{(-1)}\left(D(S,\varepsilon/2)\cap
D(0,1/\varepsilon)\right)=
\theta^{(-1)}\left(D(S,\varepsilon/2)\right) \cap
\theta^{(-1)}\left(D(0,1/\varepsilon)\right)=
\Gamma_\varepsilon\cap\Delta_\varepsilon
\]
and hence
\begin{equation} \label{2.5}
1-\limsup_{n\to\infty}\frac{q_\varepsilon(n,S)}{n}\geq
\frac{m(\Gamma_\varepsilon\cap\Delta_\varepsilon)}{m(G)}\,.
\end{equation}
By the definition of the essential range the right-hand side in
(\ref{2.5}) tends to 1 as $\varepsilon\to 0$, and so
$\lim_{n\to\infty} n^{-1}q_\varepsilon(n,S)=0$, as needed.
\\
{\bf b)}. Let $s\in S$ and $\varepsilon>0$. Construct
$F_\varepsilon$ by
\[
F_\varepsilon(z)=
  \begin{cases}
    1, & \text{for} \ z\in D(s,\varepsilon), \\
    0, & \text{for} \ z\in \mathbb C\backslash
    (D(s,2\varepsilon)),
    \end{cases} \qquad 0\leq F_\varepsilon\leq 1.
\]
Since $F_\varepsilon$ is dominated by the characteristic function
of $D(s,2\varepsilon)$, we see that
\[
\frac{\#\{\lambda\in \Sigma_n:\ \lambda\in
D(s,2\varepsilon)\}}{n}\ge \Sigma_{\lambda}(F_\varepsilon,A_n).
\]
But $\{A_n\}\sim_{\lambda} (\theta,G)$, and so by employing
$F_\varepsilon$ as the test function we have
\[
\lim_{n\rightarrow \infty}\Sigma_{\lambda}(F_\varepsilon,A_n)=
\frac1{m(G)}\int_G F_\varepsilon(\theta(t))\,dt \ge
\frac{m\{\theta^{(-1)}(D(s,\varepsilon))}{m(G)}\,,
\]
since $F_\varepsilon$ dominates the characteristic function of
$D(s,\varepsilon)$. By the definition of the essential range the
right-hand side is strictly positive and hence
\[
\liminf_{n\to\infty}\frac{\#\{\lambda\in \Sigma_n:\ \lambda\in
D(s,2\varepsilon)\}}{n}>0.
\]
The latter means exactly that $s$ attracts $\Sigma_n$ with order
$r(s)=\infty$, as was to be proved. \hfill $\square$

The final result of this Section demonstrates the stability of the
clustering under certain perturbations (cf. \cite [Corollary
4.1]{taud2}).
\begin{proposition}\label{herm-pert}
Let $\{X_n\}$ and $\{Y_n\}$ be two Hermitian matrix sequences, at
least definitely, $M$ a closed subset of the real line, and assume
that $\|X_n-Y_n\|_1=o(n)$ \ $\bigl(\|X_n-Y_n\|_1=O(1)\bigr)$. Then
$\{X_n\}$ is weakly $($strongly$)$ clustered at $M$ if and only if so
is $\{Y_n\}$.
\end{proposition}
\noindent {\em Proof}. Let $\lambda_j(X_n)$, $\lambda_j(Y_n)$ be
the eigenvalues of $X_n$ and $Y_n$, respectively, labelled in the
decreasing order. For an arbitrary $\varepsilon>0$ we introduce
three sets of indices
\begin{eqnarray*}
  I(X_n,\varepsilon) &=& \{j=1,2,\ldots,n: {\rm dist}(\lambda_j(X_n), M)>
  \varepsilon\},  \\
  I(Y_n,\varepsilon) &=& \{j=1,2,\ldots,n: {\rm dist}(\lambda_j(Y_n), M)>
  \varepsilon\}, \\
  I(X_n,Y_n,\varepsilon) &=& \{j=1,2,\ldots,n:
  |\lambda_j(X_n)-\lambda_j(Y_n)|>\varepsilon\}.
\end{eqnarray*}
Denote by $|I(X_n,\varepsilon)|$, $|I(Y_n,\varepsilon)|$ and
$|I(X_n,Y_n,\varepsilon)|$ their cardinalities. It is clear that
\[
I(X_n,\varepsilon)\subset
I\left(X_n,Y_n,\frac{\varepsilon}2\right)\bigcup
I\left(Y_n,\frac{\varepsilon}2\right),
\]
and so
\[
\left|I(X_n,\varepsilon)\right|\leq
\left|I\left(X_n,Y_n,\frac{\varepsilon}2\right)\right| +
\left|I\left(Y_n,\frac{\varepsilon}2\right)\right|.
\]
According to the Lidskii--Mirsky--Wielandt Theorem (see \cite
[Theorem IV.3.4 and Example IV.3.5] {bhatia})
\[
\sum_{j=1}^n\left|\lambda_j(X_n)-\lambda_j(Y_n)\right|\leq
\|X_n-Y_n\|_1,
\]
so that
\[
\varepsilon\,\left|I\left(X_n,Y_n,\varepsilon\right)\right|<
\sum_{j\in
I(X_n,Y_n,\varepsilon)}\left|\lambda_j(X_n)-\lambda_j(Y_n)\right|\leq
\|X_n-Y_n\|_1.
\]
Hence
\[
\left|I(X_n,\varepsilon)\right|\leq
\frac2{\varepsilon}\,\|X_n-Y_n\|_1 +
\left|I\left(Y_n,\frac{\varepsilon}2\right)\right|.
\]
The rest is plain. \hfill $\square$

\section{Non Hermitian perturbations of Hermitian matrix sequences}
\label{sec:results}

First we recall the definition of real and imaginary parts of a
matrix. Given a square matrix $A$, we define $\Ree(A)$ and
$\Imm(A)$ as $(A+A^*)/2$ and $(A-A^*)/(2i)$, respectively, where
$X^*$ denotes the conjugate transpose of the matrix $X$. In this
way, in analogy to the complex field, we naturally have
$A=\Ree(A)+i\,\Imm(A)$.

The result below is the well-known Ky Fan--Mirski Theorem (see,
e.g., \cite [Proposition III.5.3] {bhatia}).
\begin{theorem}\label{mirsky}
Let $A$ be a square matrix of size $n$, and $\lambda_j(A)$,
$\lambda_j(\Imm(A))$ the eigenvalues of $A$ and $\Imm(A)$,
respectively, labelled in the decreasing order, so that
$\Imm(\lambda_1(A))\ge \Imm(\lambda_2(A)) \ge \cdots \ge
\Imm(\lambda_n(A))$ and $\lambda_{1}(\Imm(A))\ge
\lambda_2(\Imm(A)) \ge \cdots \ge \lambda_n(\Imm(A))$. Then
\begin{equation}\label{mirsky-ineq}
\sum_{j=1}^q\Imm(\lambda_j(A))\le \sum_{j=1}^q
\lambda_j(\Imm(A)),\qquad  q=1,\ldots,n,
\end{equation}
and the equality prevails for $q=n$.\\
Equivalently, let $\lambda_j(A)$ and $\lambda_j(\Ree(A))$,
$j=1,\ldots,n$, be the eigenvalues of $A$ and $\Ree(A)$,
respectively, labelled in the decreasing order, so that
$\Ree(\lambda_1(A))\ge \Ree(\lambda_2(A)) \ge \cdots \ge
\Ree(\lambda_n(A))$ and $\lambda_{1}(\Ree(A))\ge
\lambda_2(\Ree(A)) \ge \cdots \ge \lambda_n(\Ree(A))$. Then
\begin{equation}\label{mirsky-ineq-bis}
\sum_{j=1}^q\Ree(\lambda_j(A))\le \sum_{j=1}^q
\lambda_j(\Ree(A)),\qquad q=1,\ldots,n,
\end{equation}
and the equality prevails for $q=n$.\\
\end{theorem}
The next statement provides a simple bound for the number of
nonreal eigenvalues of a matrix $A$. In what follows $\Sigma(X)$
always stands for the set of all eigenvalues of a matrix $X$:
$\Sigma(X)=\{\lambda_j(X)\}_{j=1}^n$.

\begin{lemma}\label{nonreal}
Let $A=\Ree(A)+i\,\Imm(A)$. Then for an arbitrary $\varepsilon>0$
\begin{equation}\label{3.1}
\#\{\lambda\in\Sigma(A): |\Imm (\lambda)|>\varepsilon\}\le
\frac{\|\Imm(A)\|_1}{\varepsilon}\,. \end{equation} Moreover, if
for some real $c$, $d$ we have $c\le\lambda_j(\Ree(A))\le d$ for
all $j$, then $c\le\Ree(\lambda_j(A))\le d$ and
\begin{equation}\label{3.2}
\#\{\lambda\in\Sigma(A): \lambda\notin D([c,d],\varepsilon)\}\le
\frac{\|\Imm(A)\|_1}{\varepsilon}\,.
\end{equation}
\end{lemma}
\noindent {\em Proof}. Denote by
\[
m^+:=\sum_{\lambda\in\Sigma(\Imm(A)),\ \lambda\ge 0}\lambda,
\qquad \left( m^-:=\sum_{\lambda\in\Sigma(\Imm(A)),\
\lambda<0}|\lambda|\right)
\]
the positive (negative) mass of the eigenvalues of $\Imm(A)$.
Since $\Imm(A)$ is Hermitian, its trace-norm equals the sum of the
absolute values of its eigenvalues, so $\|\Imm(A)\|_1=m^++m^-$. We
apply the first part of Theorem \ref{mirsky} for $A$ and $-A$ to
obtain
\begin{equation}\label{ineq1}
r^+:=\sum_{\lambda\in \Sigma(A),\ \Imm(\lambda)\ge 0}
\Imm(\lambda) \le m^+, \qquad r^-:=\sum_{\lambda\in \Sigma(A),\
\Imm(\lambda)< 0} |\Imm(\lambda)| \le m^-.
\end{equation}
Therefore, if we take an arbitrary $\varepsilon>0$, the number of
the eigenvalues of $A$ whose imaginary part is bigger that
$\varepsilon$ has to be bounded by $\|\Imm(A)\|_1/\varepsilon$.
Indeed,
\begin{eqnarray*}
\|\Imm(A)\|_1 &=& m^+ + m^- \ge r^+ + r^- =\sum_{\lambda\in
\Sigma(A)} |\Imm(\lambda)|\ge\sum_{\lambda\in \Sigma(A),\
|\Imm(\lambda)|>\varepsilon}|\Imm(\lambda)| \\
&\ge &\sum_{\lambda\in \Sigma(A),\ |\Imm(\lambda)|>\varepsilon}
\varepsilon
 = \varepsilon \cdot\#\{\lambda\in \Sigma(A),\ |\Imm(\lambda)|>\varepsilon\},
\end{eqnarray*}
as needed.

Next, let $\lambda$ be an eigenvalue of $A$ corresponding to an
eigenvector $\bf x$. Then
\[
\lambda={{\bf x}^*A{\bf x}\over {\bf x}^*{\bf x}}= {{\bf
x}^*\Ree(A){\bf x}\over {\bf x}^*{\bf x}}+ i{{\bf x}^*\Imm(A){\bf
x}\over {\bf x}^*{\bf x}}
\]
which implies that $\Ree(\lambda)\in [c,d]$, since, by the
assumption, every eigenvalue of $\Ree(A)$ belongs to $[c,d]$. So
(\ref{3.2}) follows from (\ref{3.1}). \hfill $\square$

\begin{corollary}\label{cor1}
Let $\{A_n\}$ be a matrix sequence such that
$\|\Imm(A_n)\|_1=o(n)$ as $n\to\infty$. Then
$q_{\varepsilon}(n,\mathbb R)=o(n)$, so $\{A_n\}$ is weakly
clustered at $\mathbb R$. Moreover, if all the eigenvalues of
$\Ree(A_n)$ are in $[c,d]$, then all the eigenvalues of $A_n$ have
real parts in the same interval and
$q_{\varepsilon}(n,[c,d])=o(n)$. The same result holds if $o(n)$
is replaced by $O(1)$ and ``weakly clustered'' by ``strongly
clustered''.
\end{corollary}

The following result establishes a link between distributions of
the Hermitian sequence $\{\Ree(A_n)\}$ and the sequence $\{A_n\}$.
As a matter of fact, we will prove a more general statement
concerning non-Hermitian perturbations of Hermitian matrix
sequences. As usual, $\|X\|$ stands for the operator (spectral)
norm of a matrix $X$.
\begin{theorem}\label{mirsky-cons}
Let $\{B_n\}$ and $\{C_n\}$ be two matrix sequences, $B_n$ is
Hermitian, at least definitely, and $A_n=B_n+C_n$. Assume further
that $\{B_n\}$ is distributed as $(\theta,G)$, $G$ of finite and
positive Lebesgue measure, both $\|B_n\|$ and $\|C_n\|$ are
uniformly bounded by a positive constant $C$ independent of $n$,
and $\|C_n\|_1=o(n)$, $n\to\infty$. Then $\theta$ is real valued
and $\{A_n\}$ is distributed as $(\theta,G)$ in the sense of the
eigenvalues. In particular, if $S(\theta)$ is the essential range
of $\theta$, then $\{A_n\}$ is weakly clustered at $S(\theta)$,
and $S(\theta)$ strongly attracts the spectra of $\{A_n\}$ with
infinite order of attraction for any of its points.
\end{theorem}
\noindent {\it Proof.} Denote by ${\rm tr} X$ the trace of a
matrix $X$, that is, the sum of its diagonal entries (or the sum
of its eigenvalues)
\[
{\rm tr} X=\sum_{\lambda\in\Sigma(X)} \lambda=\sum_{k=1}^n
(X)_{k,k},
\]
so ${\rm tr} A_n-{\rm tr} B_n={\rm tr} C_n$. As $\left|{\rm tr}
X\right|\le \|X\|_1$, the assumption on the trace-norm of $C_n$
yields
\[
\frac1{n}\,\sum_{\lambda\in \Sigma(A_n)} \lambda =\frac1{n}\,
\sum_{\lambda\in \Sigma(B_n)} \lambda + o(1).
\]
The latter is closely related to (\ref{distribution:sv-eig}) with
$F(z)=z$ (defined over the whole $\mathbb C$). Since $\{B_n\}$ is
distributed as $\theta$ over $G$, we infer by
(\ref{distribution:sv-eig})
\begin{equation}\label{from sn to an}
\lim_{n \to \infty} \frac{1}{n} \sum_{\lambda\in \Sigma(A_n)}
\lambda = \lim_{n \to \infty} \frac{1}{n} \sum_{\lambda\in
\Sigma(B_n)} \lambda
     = \frac1{m(G)}\int_G F(\theta(t)) \, dt,\ \ \ F(z)=z,
\end{equation}
where we are allowed to take $F(z)=z$ (which has an unbounded
support), since by the premises of the Theorem $\|A_n\|\le 2C$ for
all $n$, and so the spectra of $\{A_n\}$, $\{B_n\}$, and $\{C_n\}$
are all contained in the closed disk $\{|z|\le 2C\}$. Relation
(\ref{from sn to an}) can be viewed as the first step from a
distribution relation for the Hermitian sequence $\{B_n\}$ to the
same distribution relation for the sequence $\{A_n\}$. The next
step is to extend (\ref{from sn to an}) to the case when $F$ is an
arbitrary polynomial of a fixed degree. By the linearity it
suffices to consider only monomials. Clearly, for any fixed
nonnegative integer $q$, the matrix $A_n^q$ can be written as
$A_n^q=B_n^q+R_{n,q}$ and, thanks to the H\"older type
inequalities for the Schatten $p$ norms $\|X Y\|_1\le \|X\|\cdot
\|Y\|_1$ (see \cite [Corollary IV.2.6] {bhatia}) we have
$\|R_{n,q}\|_1=o(n)$ as $n\to\infty$. Therefore, by repeating the
same reasoning as above we deduce
\begin{equation}\label{from sn to an bis}
\lim_{n \to \infty} \frac{1}{n} \sum_{\lambda\in \Sigma(A_n)}
\lambda^q = \lim_{n \to \infty} \frac{1}{n} \sum_{\lambda\in
\Sigma(B_n)} \lambda^q
     = \frac1{m(G)}\int_G F(\theta(t)) \, dt,\qquad F(z)=z^q.
\end{equation}

To go over in (\ref{from sn to an bis}) from polynomials to
arbitrary continuous functions with bounded support we would like
to invoke Theorem \ref{mergelyan-cons}. So let us make sure that
the rest of its hypothesis is satisfied. As we have already
mentioned, $\|A_n\|\le 2C$ for all $n$. Next, it is clear that
\begin{equation}\label{basic-rel}
\|\Ree(C_n)\|_1\le \|C_n\|_1=o(n),\qquad \|\Imm(C_n)\|_1\le
\|C_n\|_1=o(n)
\end{equation}
as $n\to\infty$. Write $A_n=B_n+\Ree(C_n)+i\,\Imm(C_n)$. By
Theorem \ref{th-relations-def} $\{B_n\}$ is weakly clustered at
$S(\theta)$, and so is $\{\Ree(A_n)=B_n+\Ree(C_n)\}$ by Proposition
\ref{herm-pert}. Note that $S(\theta)$ is now a compact set which
lies in the interval $[-2C,2C]$, and all the eigenvalues of
$\Ree(A_n)$ are in the same interval. Corollary \ref{cor1} now
claims that $\{A_n\}$ is weakly clustered at $[-2C,2C]\supset
S(\theta)$, and the application of Theorem \ref{mergelyan-cons}
completes the proof. \hfill $\square$

The following theorem deals with the case of the strong
clustering.

\begin{theorem}\label{mirsky-cons-bis}
Let $\{B_n\}$ and $\{C_n\}$ be two matrix sequences, $B_n$ is
Hermitian, at least definitely, and $A_n=B_n+C_n$. Assume that
$\{B_n\}$ is strongly clustered at $[c,d]$, $\|C_n\|_1=O(1)$,
$n\to\infty$ and $\|A_n\|$ is uniformly bounded by a positive
constant $C$ independent of $n$. Then $\{A_n\}$ is strongly
clustered at $[c,d]$.
\end{theorem}
\noindent {\it Proof.} Since now $\|\Ree(C_n)\|_1=O(1)$ and
$\|\Imm(C_n)\|_1=O(1)$, both the related sequences are strongly
clustered at zero by Proposition \ref{herm-pert}. A repeated
application of the same proposition shows that $\{B_n+\Ree(C_n)\}$
is strongly clustered at $[c,d]$. Although we don't have the right
to invoke Corollary \ref{cor1} at this point, since the
eigenvalues of $\Ree(A_n)=B_n+\Ree(C_n)$ are not necessarily in
$[c,d]$, we can follow a direct approach stemming from Theorem \ref{mirsky}.

Since $\|A_n\|\le C$, the real part of any eigenvalue of $A_n$
belongs to $[-C,C]$ and the same is true for any eigenvalue of
$\Ree(A_n)$. For $\varepsilon>0$, let $q_n^-(\varepsilon)$ be the
number of eigenvalues of $A_n$ whose real parts are below
$c-\varepsilon$, and analogously, let $q_n^+(\varepsilon)$ be the
number of eigenvalues of $X_n$ whose real parts exceed
$d+\varepsilon$. We want to prove that both $q_n^-(\varepsilon)$
and $q_n^+(\varepsilon)$ can be bounded by a constant possibly
depending on $\varepsilon$, but independent of $n$. By
(\ref{mirsky-ineq-bis}) we have
\[
\sum_{j=1}^{q_n^+(\varepsilon)}\Ree(\lambda_j(A_n))\le
\sum_{j=1}^{q_n^+(\varepsilon)} \lambda_j(\Ree(A_n))
\]
with
\[
\Ree(\lambda_1(A_n))\ge \Ree(\lambda_2(A_n)) \ge \cdots \ge
\Ree(\lambda_{q_n^+(\varepsilon)}(A_n))>d+\varepsilon\ge
\Ree(\lambda_{q_n^+(\varepsilon)+1}(A_n)).
\]
Therefore
\begin{equation}\label{mirsky-ineq-conseq}
(d+\varepsilon)q_n^+(\varepsilon) \le
\sum_{j=1}^{q_n^+(\varepsilon)} \lambda_j(\Ree(A_n)).
\end{equation}
Thanks to the strong clustering of $\Ree(A_n)=B_n+\Ree(C_n)$, for
every $\varepsilon'>0$ there exists a positive constant
$K(\varepsilon')$ independent of $n$ such that the number of
eigenvalues of $\Ree(A_n)$ not belonging to
$(c-\varepsilon',d+\varepsilon')$ is bounded by $K(\varepsilon')$.
Consequently, we infer
\begin{equation}\label{mirsky-ineq-conseq-bis}
\sum_{j=1}^{q_n^+(\varepsilon)} \lambda_j(\Ree(A_n)) \le
CK(\varepsilon')+
(d+\varepsilon')(q_n^+(\varepsilon)-K(\varepsilon'))^+
\end{equation}
with $(x)^+=(x+|x|)/2$. Putting together
(\ref{mirsky-ineq-conseq}) and (\ref{mirsky-ineq-conseq-bis}), by
choosing $\varepsilon'=\varepsilon/2$, we finally deduce
\[
q_n^+(\varepsilon)\le {2CK(\varepsilon/2)\over \varepsilon}\,,
\]
where, as requested, the right-hand side is independent of $n$. A
similar reasoning on $-X_n$ gives the same bound on
$q_n^-(\varepsilon)$, as claimed.

As for the imaginary parts of the eigenvalues of $A_n$, we can
apply directly (\ref{3.1}). The proof is complete. \hfill
$\square$

The latter result can be extended to the case of clustering at
several intervals, the situation we will encounter later in
Theorem \ref{wcesaro-bis}.
\begin{theorem}\label{sever-inter}
Let $\{B_n\}$ and $\{C_n\}$ be two matrix sequences, $B_n$ is
Hermitian, at least definitely, and $A_n=B_n+C_n$. Let $E$ be a
union of $m$ disjoint closed intervals $($possibly, degenerate$)$.
Assume that $\{B_n\}$ is strongly clustered at $E$, $\|C_n\|_1=O(1)$,
$n\to\infty$ and $\|A_n\|$ is uniformly bounded by a positive
constant $C$ independent of $n$. Then $\{A_n\}$ is strongly
clustered at $E$.
\end{theorem}
\noindent {\it Proof.} We reduce this statement to the previous
one. Denote
\[
E=\cup_{j=1}^m [a_j, b_j], \qquad a_1\le b_1<a_2\le
b_2<\ldots<a_m\le b_m,
\]
and put $T(z)=\prod_{j=1}^m (z-a_j)(z-b_j)$. Obviously, $T(E)\in
[\omega,0]$, $\omega=\min_x T(x)<0$, and $T(x)>0$ for $x\in\mathbb
R\backslash E$. By the Spectral Mapping Theorem (see e.g. \cite
[p. 5]{bhatia}) $E$ is a strong cluster for $\{B_n\}$ yields
$[\omega,0]$ is a strong cluster for $\{T(B_n)\}$. Next, by the
hypothesis of the theorem and the H\"older type inequalities for
the trace norm
\[
T(A_n)=T(B_n)+R_n, \qquad \|R_n\|_1=O(1), \ n\to\infty.
\]
We have the right to apply Theorem \ref{mirsky-cons-bis} to the
matrix sequences $\{T(A_n)\}$, $\{T(B_n)\}$ to conclude that $\{T(A_n)\}$ is
strongly clustered at $[\omega,0]$. The repeated application of
the Spectral Mapping Theorem completes the proof. \hfill $\square$

Let us go back to the Jacobi matrix sequences described in the
introduction. Let
\[
A_\infty = \left[ \begin{array}{ccccc}
    b_0 & c_1 & \\
    a_1 & b_1 & c_2 & \\
        & a_2 & b_2& c_3 \\
        &   &\ddots & \ddots & \ddots
             \end{array} \right]
\]
be an infinite complex non symmetric Jacobi matrix with the
bounded entires
\begin{equation}\label{3.10}
\sup_n\,(|a_n|+|b_n|+|c_n|)\le C<\infty.
\end{equation}
As a simple consequence of Theorem \ref{mirsky-cons}, we can prove
the following

\begin{corollary}\label{cor2}
Let $A_\infty$ be the Ces\`aro compact perturbation of
$J_\infty^0$, that is,
\begin{equation}\label{cesaro}
\lim_{n\to\infty} \frac1n\,\sum_{j=1}^n(|1-a_j|+|b_j|+|1-c_j|)=0,
\end{equation}
and $\{A_n\}$ its principal $n\times n$ blocks. Then $\{A_n\}$ is
distributed as $(2\cos t, [-\pi,\pi])$ in the sense of
eigenvalues, weakly clustered at $[-2,2]$, and $[-2,2]$ strongly
attracts the spectra of $\{A_n\}$ with infinite order of
attraction for any of its points.
\end{corollary}
\noindent {\em Proof.} We apply Theorem \ref{mirsky-cons} with
$B_n=J_n^0$, $C_n=A_n-J_n^0$. $\{B_n\}$ is clearly distributed as
$(2\cos t, [-\pi,\pi])$. Inequality (\ref{3.10}) provides the
uniform boundedness of $\|A_n\|$ and $\|C_n\|$. Finally,
(\ref{cesaro}) is equivalent to $\|C_n\|_1=o(n)$, and the result
follows. \hfill $\square$

\begin{corollary}\label{cor2-bis}
Let $A_\infty$ be trace class perturbation of $J_\infty$,
that is,
\begin{equation}\label{trace-class}
\limsup_{n\to\infty} \sum_{j=1}^n(|1-a_j|+|b_j|+|1-c_j|)<\infty.
\end{equation}
Then $\{A_n\}$ is distributed as $(2\cos t, [-\pi,\pi])$ in the
sense of eigenvalues, strongly clustered at $[-2,2]$, and $[-2,2]$
strongly attracts the spectra of $\{A_n\}$ with infinite order of
attraction for any of its points.
\end{corollary}
\noindent {\em Proof.} The only point to be proved is that the
weak cluster is also strong, and this is implied by Theorem
\ref{mirsky-cons-bis}. \hfill $\square$

If we are concerned only about the clustering of the spectrum, the
more elementary Corollary \ref{cor1} does the job. In this case
the assumption $\Imm(A_\infty)$ being a Ces\`aro compact
perturbation of $J_\infty^0$, that is,
\[
\lim_{n\to\infty} \frac1n\,\sum_{j=1}^n(|\Imm( b_j)|+|a_j-\bar
c_j|)=0,
\]
already guarantees that $\{A_n\}$ is weakly clustered at $\mathbb
R$. Moreover, if all the eigenvalues of $\Ree(A_n)$ are in
$[c,d]$, then $\{A_n\}$ is weakly clustered at $[c,d]$.

It is worth pointing out that $\Sigma(A_n)$ now agrees with the
set of all zeros of the polynomial $p_n$ which satisfies the
three-term recurrence relation
\begin{equation}\label{ops}
zp_j(z)=a_{j}p_{j-1}(z)+b_jp_j(z)+c_{j+1}p_{j+1}(z), \qquad
j\in\mathbb Z_+
\end{equation}
$p_{-1}=0$, $p_0=1$. Such polynomials are studied systematically
in the theory of Pad\'e approximations and continued
$J$-fractions. More precisely, $p_n$ is the denominator of the
$n$th diagonal Pad\'e approximant and its zeros are the poles of
this Pad\'e approximant. In turn, the closed interval $[-2,2]$ is
now the essential spectrum of the bounded operator $A_\infty$ in
$\ell^2$.

{\sl Remark}. In \cite{blmt1, blmt2} the authors studied the
attracting properties of the spectrum of $A_\infty$ in the case
when $A_\infty$ is a compact perturbation of $J_\infty$. The
celebrated theorem of H. Weyl claims that
$\Sigma(A_\infty)=[-2,2]\cup\Sigma_d(A_\infty)$, where the
discrete spectrum $\Sigma_d(A_\infty)$ is at most denumerable set
of eigenvalues $\lambda_j(A_\infty)$ of the finite algebraic
multiplicity $\nu_j$, off the essential spectrum $[-2,2]$. It is
proved in \cite{blmt1, blmt2} that each $\lambda_j$ is the
attracting point of $\Sigma(A_n)$ of order $\nu_j$. Our result in
Corollary \ref{cor2} supplements this one nicely. Note that in the
case (\ref{cesaro}) the Weyl theorem is only partly true (see
\cite [Theorems 7 and 9]{golin}): still
$[-2,2]\subset\Sigma(A_\infty)$, but in general there is no
discrete part of the spectrum any more.

\section{Asymptotically periodic Jacobi matrices: the block
case}\label{sec:periodic}

We start out with the definition of the spectral distribution with
matrix-valued symbols. Throughout the rest of the paper $\theta$
will stand for a $k\times k$ matrix-valued and the Lebesgue
integrable function (i.e., all its entries are integrable) with
the eigenvalues $\lambda_j(\theta)$, $j=1,2,\ldots,k$.

\begin{definition}   \label{def:asd}
{\rm Let $\theta$ be a $k\times k$ matrix-valued the Lebesgue
integrable function defined on a set $G$ of finite Lebesgue
measure. A matrix sequence $\{ A_{n} \}$ has the {\em asymptotic
spectral distribution} $\theta$ if for all $F \in \mathcal{C}_{0}$
one has
\[
   \lim_{n\to \infty}
   \Sigma_{\lambda}(F,A_n) =
   \frac1{km(G)} \sum_{j=1}^k
   \int_{G}F(\lambda_j(\theta(t)))\,dt.
\]
As in the scalar case, we write in short} $\{ A_{n} \}\sim_\lambda
(\theta,G)$.
\end{definition}
Under the essential range of $\theta$ we mean now the set
\[
S(\theta):= \bigcup_{j=1}^k {\rm Range}\,(\lambda_j(\theta)).
\]
The same argument as applied above in the proof of Theorem
\ref{th-relations-def} leads to the following result.

\begin{theorem}\label{th-relations-def1}
Let $\theta$ be a $k\times k$ matrix-valued, the Lebesgue
integrable function defined on a set $G$ of finite Lebesgue
measure and $S=S(\theta)$ be the essential range of $\theta$. Let
$\{A_n\}$ be a matrix sequence distributed as $\theta$ in the
sense of eigenvalues. Then
\begin{description}
\item[a)]  $S(\theta)$ is a weak cluster for $\{A_n\}$; \item[b)]
each point $s\in S(\theta)$ strongly attracts $\Sigma_n$ with
infinite order $r(s)=\infty$.
\end{description}
\end{theorem}

\begin{definition}\label{toeplitz-block}
{\rm Let $b$ be a $k \times k$ matrix-valued and the Lebesgue
integrable function defined on $[-\pi,\pi)$ with the Fourier
coefficients
\begin{equation}   \label{eq:fourier}
    {\hat b}_j =\frac1{2\pi}
    \int_{-\pi}^\pi b(t) \exp(- ijt)\,dt
    \in \mathbb{C}^{k \times k}, \qquad j \in \mathbb Z.
\end{equation}
The function $b$ is called the {\em generating function} of the
sequence of block Toeplitz matrices}
$$
   T_{n}(b) = \left[ \begin{array}{c c c c}
   {\hat b}_0 & {\hat b}_{-1} & \cdots & {\hat b}_{1-n} \\
   {\hat b}_{1} & {\hat b}_0 & \ddots & \vdots \\
   \vdots & \ddots & \ddots & {\hat b}_{-1} \\
   {\hat b}_{n-1} & \cdots & {\hat b}_{1} & {\hat b}_0
   \end{array} \right]\in\mathbb C^{kn\times kn}.
$$
\end{definition}

It is easy to observe that $T_{n}(b)$ is Hermitian for every $n$
if and only if its generating function $b$ is Hermitian for almost
every $t\in [-\pi,\pi)$, and the index $n$ here denotes the block
order.

Let us define a matrix sequence $\{\tilde T_m(b)\}$ by the
following recipe: $\tilde T_{kn}:=T_n$, and $\tilde T_{kn-j}$ is
obtained from $T_n$ by deleting the last $j$ rows and columns for
$j=1,2,\ldots,k-1$. In other words, $\tilde T_m$ is the principal
$m\times m$ block of the infinite block-matrix $T_\infty(b)=\{\hat
b_{p-q}\}_{p,q=0}^\infty$.

The following general result due to Tilli (see \cite{Tillinota})
is very important in our context.

\begin{theorem}   \label{theor:generating}
If $b$ is any Hermitian-valued and absolutely integrable function
on $[-\pi,\pi]$
$$
   \int_{-\pi}^\pi \| b(t) \| \, dt < +\infty,
$$
where $\| \, \cdot \, \|$ is any matrix norm in $\mathbb{C}^{k
\times k}$, then $\{\tilde T_m(b)\}\sim_\lambda (b,[-\pi,\pi])$ in
the sense of Definition~$\ref{def:asd}$.
\end{theorem}

Now we turn to the case of asymptotically periodic Jacobi
matrices. Let
\begin{equation} \label{4.1}
J_\infty^{(0)} = \left[ \begin{array}{ccccc}
    b_0^{(0)} & a_1^{(0)} & \\
    a_1^{(0)} & b_1^{(0)} & a_2^{(0)} & \\
        & a_2^{(0)} & b_2^{(0)} & a_3^{(0)} \\
        &   &\ddots & \ddots & \ddots
             \end{array} \right], \qquad a_n^{(0)}>0, \quad
             b_n^{(0)}\in\mathbb R
\end{equation}
be an infinite Jacobi matrix with $k$-periodic entries
\begin{equation} \label{4.2}
a_{n+k}^{(0)}=a_n^{(0)}, \qquad  b_{n+k}^{(0)}=b_n^{(0)}, \qquad
n\in\mathbb Z_+,
\end{equation}
and ${\bf a}  = (a_0^{(0)}, a_1^{(0)}, \ldots , a_{k-1}^{(0)})$,
${\bf b} = (b_0^{(0)}, b_1^{(0)}, \ldots , b_{k-1}^{(0)})$ be two
real vectors of order $k$, which define completely the entries of
the whole matrix $J_\infty^{(0)}$. In that case the principal
$m\times m$ block of $J_\infty^{(0)}$ in (\ref{4.1}) is denoted
more explicitly by $J_m^{(0)}=J_m[{\bf a},{\bf b}]$.

Let $\theta({\bf a}, {\bf b}, t)$ be the Hermitian matrix-valued
trigonometric polynomial of the form
\begin{equation}   \label{eq:fab}
   \theta({\bf a}, {\bf b}, t) =
   J_{k}[{\bf a},{\bf b}] +
   \left[ \begin{array}{c c c c}
   0 & \cdots & 0 & a_0^{(0)}\exp(it) \\
   \vdots  & \cdots & 0 & 0 \\
   0 & 0 & \cdots & \vdots \\
   a_0^{(0)}\exp(-it) & 0 & \cdots & 0
   \end{array} \right] . 
\end{equation}
It is a matter of simple computation to verify that $\theta$ has
only three nonzero Fourier coefficients $\hat\theta_0$ and
$\hat\theta_{\pm 1}$,
\[
T_n(\theta({\bf a}, {\bf b}))=
 \left[ \begin{array}{c c c c}
   {\hat \theta}_0 & {\hat \theta}_{-1} & & O \\
   {\hat \theta}_1 & {\hat \theta}_0 &
         \ddots & \\
       & \ddots & \ddots & {\hat b}_{-1} \\
   O   & & {\hat \theta}_1 & {\hat \theta}_{0}
   \end{array} \right]=J_{kn}[{\bf a},{\bf b}],
\]
and in general $\tilde T_m(\theta({\bf a}, {\bf b}))=J_{m}[{\bf
a},{\bf b}]$ for all $m$. So the asymptotic distribution for the
periodic Jacobi matrix sequence is a particular case of Theorem
\ref{theor:generating} with $b=\theta({\bf a}, {\bf b})$. Such
asymptotic distribution, paraphrased as the asymptotic
distribution of the zeros of orthogonal polynomials $p_n$
(\ref{ops}) with periodic recurrence coefficients, is well known
(see, e.g., \cite [Section 3]{assche} and references therein). The
essential range $S(\theta({\bf a}, {\bf b}))$ is tightly related
to the support of the corresponding orthogonality measure (cf.
\cite [Theorem 13]{mnas}). If $k=1$, then $\theta({\bf a},{\bf b})
= b^{(0)} + 2a^{(0)}\cos t$ and putting $b^{(0)}=0$ and
$a^{(0)}=1$ we come to the Toeplitz matrix (\ref{toep-conv}).

As in the scalar case ($k=1$), we are interested in generic
complex perturbations of $J_\infty^{(0)}$. An infinite complex
Jacobi matrix
\begin{equation} \label{4.5}
J_\infty = \left[ \begin{array}{ccccc}
    b_0 & c_1 & \\
    a_1 & b_1 & c_2 & \\
        & a_2 & b_2 & c_3 \\
        &   &\ddots & \ddots & \ddots
             \end{array} \right], \qquad a_n, b_n, c_n\in\mathbb C
\end{equation}
is called the {\em Ces\`aro asymptotically $k$-periodic} if
\[
\lim_{n\to\infty}\frac1n \sum_{j=1}^n
(|a_j-a_j^{(0)}|+|b_j-b_j^{(0)}|+|c_j-a_j^{(0)}|)=0,
\]
the {\em asymptotically $k$-periodic} if
\[
\lim_{n\to\infty}(|a_n-a_n^{(0)}|+|b_n-b_n^{(0)}|+|c_n-a_n^{(0)}|)=0,
\]
and the {\em trace class asymptotically $k$-periodic} if
\[
\limsup_{n\to\infty}\sum_{j=1}^n(|a_j-a_j^{(0)}|+|b_j-b_j^{(0)}|+
|c_j-a_j^{(0)}|)<\infty,
\]
for some $k$-periodic sequences $\{a_n^{(0)}, b_n^{(0)}\}$ as in
(\ref{4.2}). In other words, $J_\infty=J_\infty^{(0)}+P_\infty$
with the $k$-periodic $J_\infty^{(0)}$ (\ref{4.1}) (called the
background) and the Ces\`aro compact (compact, the trace class)
perturbation $P_\infty$.

The following results can be proved in exactly the same fashion as
Corollary \ref{cor2} and Corollary \ref{cor2-bis}. In the latter
case Theorem \ref{sever-inter} comes into play. The point is that
the essential range $S(\theta({\bf a}, {\bf b}))$ is now a union
of at most $k$ disjoint closed intervals, and all the eigenvalues
of $J_n^{(0)}$ (the zeros of orthogonal polynomials $p_n^{(0)}$
(\ref{ops})), but finitely many (at most $2k$), lie in
$S(\theta({\bf a}, {\bf b}))$. So, in particular, the matrix
sequence $\{J_n^{(0)}\}$ is strongly clustered at $S(\theta({\bf
a}, {\bf b}))$.

\begin{theorem}\label{wcesaro}
Let $J_\infty$ be the Ces\`aro asymptotically $k$-periodic Jacobi
matrix with the background $J_\infty^{(0)}$ and $\theta({\bf a},
{\bf b})$ $(\ref{eq:fab})$ the generating function for
$J_\infty^{(0)}$. Then $\{J_n\}$ is distributed as $(\theta({\bf
a}, {\bf b}), [-\pi,\pi])$ in the sense of eigenvalues, weakly
clustered at $S(\theta({\bf a}, {\bf b}))$, and $S(\theta({\bf a},
{\bf b}))$ strongly attracts the spectra of $\{J_n\}$ with
infinite order of attraction for any of its points.
\end{theorem}

\begin{theorem}\label{wcesaro-bis}
Let $J_\infty$ be the trace class asymptotically $k$-periodic
Jacobi matrix with the background $J_\infty^{(0)}$ and
$\theta({\bf a}, {\bf b})$ $(\ref{eq:fab})$ the generating
function for $J_\infty^{(0)}$. Then $\{J_n\}$ is strongly
clustered at $S(\theta({\bf a}, {\bf b}))$, and $S(\theta({\bf a},
{\bf b}))$ strongly attracts the spectra of $\{J_n\}$ with
infinite order of attraction for any of its points.
\end{theorem}

\section{Concluding remarks and further generalizations}\label{sec:extension}

As a conclusion, we observe that tools from matrix theory \cite{bhatia,B2000} combined
with those from asymptotic linear algebra \cite{Tillinota,tillicomplex,taud2}
have been crucial for proving plainly results concerning
non Hermitian perturbation of Jacobi matrix sequences. A special part of them is the GLT
theory (see \cite{ser-glt,glt-vs-fourier} and references therein) which allows to treat the case of
variable coefficients under very mild restrictions
on the regularity of the coefficients (e.g. numerical approximations of
variable coefficient PDEs \cite{ser-glt} and systems of PDEs  \cite{glt-vs-fourier},
Jacobi sequences with asymptotically varying periodic
\cite{fasino} and non-periodic \cite{ku-ser} coefficients, etc.).
The interesting fact is that the tools explicitly developed here are applicable
verbatim to these cases as well, by allowing to deal with non Hermitian perturbations
under the same mild trace conditions.

\section{Appendix. Equivalence of trace-norm and entry-wise conditions}

Let $A=\left\{a_{j,k}\right\}_{j,k=1}^n$ be a complex matrix of
size $n$, let $\|\cdot\|_1$ be the trace-norm, and let
$\|\cdot\|_{[1]}$ be the componentwise $l^1$ norm:
\[
\|A\|_1=\sum_{j=1}^n\sigma_j,\qquad \|A\|_{[1]}=\sum_{j,k=1}^n
|a_{j,k}|
\]
with $\sigma_1\ge \sigma_2\ge \cdots \sigma_n\ge 0$ being the
singular values of $A$. With the notations (i) and (ii) at the end
of Section \ref{sec:intro}, we would like to prove that
$\|P_n\|_1=o(n)$ if and only if (\ref{equiv-1}) holds and
$\|P_n\|_1=O(1)$ if and only if (\ref{equiv-2}) is satisfied.
Taking into account the definition of the norm $\|\cdot\|_{[1]}$
and the tridiagonal structure of $P_n$, $n\ge 1$, condition
(\ref{equiv-1}) can be rewritten as $\|P_n\|_{[1]}=o(n)$ and,
similarly, (\ref{equiv-2}) is equivalent to $\|P_n\|_{[1]}=O(1)$.
Therefore what we would like to prove is the asymptotic
equivalence, independently of the size $n$, of the two norms
$\|\cdot\|_1$ and $\|\cdot\|_{[1]}$. Specifically, we look for two
positive constants $c$ and $C$ independent of $n$ such that
$c\|A\|_1\le \|A\|_{[1]}\le C \|A\|_1$ for every complex matrix
$A$ of size $n$. For a fixed $n$, the existence of the two
positive constants $c=c(n)$ and $C=C(n)$ is trivial thanks to the
topological equivalence of norms in any finite dimensional vector
space. The nontrivial part is to show that $c$ and $C$ can be
chosen independently of $n$. Unfortunately, the latter is in
general false as the following example shows. Take
$A=\left\{a_{j,k}\right\}_{j,k=1}^n$ with $a_{j,k}=1$, $\forall
j,k=1,\ldots,n$. Then $\sigma_1=n$, $\sigma_2=\cdots=\sigma_n=0$,
and therefore $\|A\|_1=n$ while $\|A\|_{[1]}=n^2$ so that $C(n)\ge
n$ (indeed it can be proved that the previous example is an
extremal one and indeed the best constant $C$ is exactly
$C(n)=n$).

Therefore the equivalence of the trace-norm and of the $l^1$
entry-wise norm has to exploit the fact that the involved matrices
are tridiagonal. In the subsequent steps we will use the Fourier
analysis of matrices introduced by Bhatia in \cite{B2000}. Let $A$
be a generic tridiagonal matrix of size $n$ and, for any
$m=1-n,\ldots,n-1$, let ${\cal D}_m(A)$ be the matrix which
coincides with the $m$-th diagonal of $A$, i.e., $\left\{{\cal
D}_m(A)\right\}_{j,k}=a_{j,k}$ if $j-k=m$ and $\left\{{\cal
D}_m(A)\right\}_{j,k}=0$ otherwise. Therefore
\begin{equation}\label{rel-1}
A=\sum_{m=-1}^1 {\cal D}_m(A)
\end{equation}
and, by the structure of any ${\cal D}_m(A)$, a plain check shows
that
\begin{equation}\label{rel-2}
\|{\cal D}_m(A)\|_{[1]}=\|{\cal D}_m(A)\|_{1}.
\end{equation}
Consequently, by the definition of $\|\cdot\|_{[1]}$,
(\ref{rel-1}), and (\ref{rel-2}) we have
\begin{eqnarray*}
\|A\|_1 & = & \|\sum_{m=-1}^1 {\cal D}_m(A)\|_1
\le  \sum_{m=-1}^1 \|{\cal D}_m(A)\|_1 \\
& = & \sum_{m=-1}^1 \|{\cal D}_m(A)\|_{[1]} = \|A\|_{[1]}
\end{eqnarray*}
and so $c=1$ which is independent of $n$. For the reverse
inequality we have
\begin{eqnarray*}
\|A\|_{[1]} & = & \|\sum_{m=-1}^1 {\cal D}_m(A)\|_{[1]}
= \sum_{m=-1}^1 \|{\cal D}_m(A)\|_{[1]} \\
& = & \sum_{m=-1}^1 \|{\cal D}_m(A)\|_1 \le 3\|A\|_1,
\end{eqnarray*}
where for the last inequality we use the identity (see
\cite{B2000})
\[
{\cal D}_m(A)=\frac1{2\pi}\,\int_{-\pi}^\pi D(t)AD^*(t)
{\exp}(-imt)\, dt,
\]
with $D(t)$ a diagonal unitary matrix whose $j$-th diagonal entry
equals ${\exp}(i(j-1)t)$. From the latter identity, since the
trace-norm is a unitarily invariant norm (see \cite{bhatia}), it
easily follows that $\|{\cal D}_m(A)\|_1\le \|A\|_1$. We conclude
that $C=3$ which is again a constant independent of $n$, as
desired.

As already pointed out in the introduction only the proof is new.
The result can be recovered directly from known facts: for
instance use inequalities (2.32) in \cite{kilsimon}  with $p=1$
and the (trivial) equivalence between $l^\infty$ and $l^1$ norms
for vectors of size $3$. Then one arrives to
\[
{1\over 3}\|A\|_1\le \|A\|_{[1]}\le 9 \|A\|_1
\]
for every tridiagonal matrix $A$. Note however that our constants
$c=1$ and $C=3$ are tighter and indeed $c=1$ is optimal (take $A$
the identity matrix).

Finally, it should be remarked that the similar equivalence
results can be obtained for more general patterns. Instead of
tridiagonal structures we could equally well have considered
banded structures (also in a multilevel sense, see \cite{HHS}). In
that case, the proofs are identical and the constants are $c=1$
and $C$ equals to the number of nonzero diagonals of the
considered band matrices. As long as this number is independent of
$n$, the two norms $\|\cdot\|_{[1]}$ and $\|\cdot\|_{1}$ are
asymptotically equivalent, i.e., with equivalence constants
positive and independent of the size~$n$.

\vfill
\eject
\noindent
Leonid Golinskii,\\
Mathematical Division,
Institute for Low Temperature Physics,\\
Kharkov University,\\
47 Lenin ave, Kharkov
61103, Ukraine;\\
e-mail: : golinskii@ilt.kharkov.ua.

\ \\
\\
\\
Stefano Serra-Capizzano,\\
Department of Physics and Mathematics,\\
University of ``Insubria'',\\
Via Valleggio 11, 22100 Como, Italy;\\
e-mail: stefano.serrac@uninsubria.it.

\end{document}